\documentclass[11pt,oneside,a4paper]{article}
\usepackage[latin1]{inputenc}
\usepackage{ulem}
\usepackage{color}
\usepackage{fancybox}
\usepackage{fancyhdr}
\usepackage{amssymb}
\usepackage{amsmath}
\usepackage{amsthm}
\usepackage{stmaryrd}
\usepackage{graphicx}
\usepackage{hyperref}
\usepackage{epic}
\usepackage{fixmath}
\usepackage{bbm}
\usepackage{supertabular}
\usepackage{longtable}
\usepackage[left=3cm,right=3cm,top=2.5cm,bottom=2cm]{geometry}

\newcommand{\HH}{\mathbb{H}}
\newcommand{\C}{\mathbb{C}}
\newcommand{\R}{\mathbb{R}}

\newcommand{\GL}{\mathrm{GL}}

\newcommand{\Sp}{\mathrm{Sp}}
\newcommand{\U}{\mathrm{U}}
\newcommand{\SU}{\mathrm{SU}}
\newcommand{\SO}{\mathrm{SO}}

\newcommand{\Spin}{\mathrm{Spin}}

\newcommand{\Phirm}{\mathrm{\Phi}}
\newcommand{\Iso}{\mathrm{Iso}}

\newcommand{\g}{\mathfrak{g}}
\newcommand{\Ad}{\mathrm{Ad}}
\newcommand{\so}{\mathfrak{so}}
\newcommand{\su}{\mathfrak{su}}
\newcommand{\ssp}{\mathfrak{sp}}
\newcommand{\cc}{\mathfrak{c}}
\newcommand{\hh}{\mathfrak{h}}
\newcommand{\uu}{\mathfrak{u}}
\newcommand{\s}{\mathfrak{s}}

\newtheorem{theorem}{Theorem}[section]
\newtheorem{proposition}{Proposition}[section]
\newtheorem{lemma}{Lemma}[section]

\newtheorem{remark}{Remark}[section]

\font\tenrm=cmr10
\font\cmssl=cmss10 at 12 pt
\font\bigss=cmssdc10 scaled 2300

\begin{document}
\rightline{}
\vskip 1.5 true cm
\begin{center}
{\bigss Pseudo-Riemannian almost quaternionic\\[.5em] homogeneous spaces with irreducible isotropy}
\vskip 1.0 true cm
{\cmssl  V.\ Cort\'es$^1$ and B.\ Meinke$^2$
} \\[3pt]
{\tenrm $^1$Department of Mathematics\\
and Center for Mathematical Physics\\
University of Hamburg\\
Bundesstra{\ss}e 55,
D-20146 Hamburg, Germany\\
vicente.cortes@uni-hamburg.de 
\\[1em]
$^2$Department of Mathematics\\ 
University of D\"usseldorf\\
Universit\"atsstra{\ss}e 1, Raum 25.22.03.62 \\
D-40225 D\"usseldorf, Germany\\  
Benedict.Meinke@uni-duesseldorf.de}
\end{center}
\vskip 1.0 true cm
\baselineskip=18pt
\begin{abstract}
\noindent
We show that pseudo-Riemannian almost quaternionic homogeneous spaces with index $4$ and an $\HH$-irreducible isotropy group are locally isometric to a pseudo-Riemannian quaternionic K\"ahler symmetric space if the dimension is at least $16$. In dimension $12$ we give a non-symmetric example.\\
{\it Keywords: Homogeneous spaces, symmetric spaces, pseudo-Riemannian manifolds, almost quaternionic structures}\\[.5em]
{\it MSC classification: 53C26, 53C30, 53C35, 53C50}
\end{abstract}

\section{Introduction}
In \cite{AhmedZeghib} Ahmed and Zeghib studied pseudo-Riemannian almost complex homogeneous spaces of index $2$ with a $\C$-irreducible isotropy group. They showed that these spaces are already pseudo-K\"ahler if the dimension is at least $8$. If furthermore the Lie algebra of the isotropy group is $\C$-irreducible then the space is locally isometric to one of five symmetric spaces.\\
There are two different quaternionic analogues of K\"ahler manifolds, namely hyper-K\"ahler and quaternionic K\"ahler manifolds. In the first case, the complex structure is replaced by three complex structures assembling into a hyper-complex structure $(I,J,K)$, in the second by the more general notion of a quaternionic structure $Q\subset \mathrm{End}\, TM$ on the underlying manifold $M$. Riemannian as well as pseudo-Riemannian quaternionic K\"ahler manifolds are Einstein and therefore of particular interest in pseudo-Riemannian geometry.\\
In \cite{CortesMeinke} the authors investigated the hyper-complex analogue of the topic studied by Ahmed and Zeghib, namely pseudo-Riemannian almost hyper-complex homogeneous spaces of index $4$ with an $\HH$-irreducible isotropy group. It turned out that these spaces of dimension greater or equal than $8$ are already locally isometric to the flat space $\HH^{1,n}$ except in dimension $12$, where non-symmetric examples exist.\\
In this article we study the quaternionic analogue, that is we consider pseudo-Rie\-mannian almost quaternionic homogeneous spaces of index $4$ with an $\HH$-irreducible isotropy group. The main result of our analysis is the following theorem. 
\begin{theorem}\label{th:MainTheorem}
  Let $(M,g,Q)$ be a connected almost quaternionic pseudo-Hermitian manifold of index $4$ and $\dim M=4n+4\geq 16$, such that there exists a connected Lie subgroup $G\subset\Iso(M,g,Q)$ acting transitively on $M$. If the isotropy group $H:=G_p$, $p\in M$, acts $\HH$-irreducibly, then $(M,g,Q)$ is locally isometric to a quaternionic K\"ahler symmetric space.
\end{theorem}
\noindent Here $\Iso(M,g,Q)$ denotes the subgroup of the isometry group $\Iso(M,g)$ which preserves the almost quaternionic structure $Q$ of $M$. A consequence of the theorem is that the homogeneous space $M$ itself is quaternionic K\"ahler and locally symmetric. Notice that pseudo-Riemannian quaternionic K\"ahler symmetric spaces have been classified in \cite{AlekseevskyCortes}. In Section \ref{counterSec} we show, by construction of a non-symmetric example in dimension $12$, that the hypothesis $\dim M \ge 16$ in Theorem \ref{th:MainTheorem} cannot be omitted. Moreover, we classify in Proposition \ref{classProp} all examples with the same isotropy algebra $\hh=\so(1,2)\oplus\so(3)\subset \so(1,2)\oplus\so(4)\subset \mathfrak{gl}(\mathbb{R}^{1,2}\otimes \mathbb{R}^4)\cong \mathfrak{gl}(12,\mathbb{R})$ in terms of the solutions of a system of four quadratic equations for six real variables.\\[3pt]
The strategy of the proof of Theorem \ref{th:MainTheorem} is as follows. We consider the $\HH$-irreducible isotropy group $H$ as a subgroup of $\Sp(1,n)\Sp(1)$ and classify the possible Lie algebras. Then we consider the covering $G/H^0$ of $M=G/H$ and show by taking into account the possible Lie algebras that it is a reductive homogeneous space. Finally, we show that the universal covering $\tilde{M}$ is a symmetric space. The invariance of the fundamental $4$-form under $G$ then implies that the symmetric space is quaternionic K\"ahler.

\noindent\textbf{Acknowledgments.} This work was partly supported by the German Science Foundation (DFG) under the
Collaborative Research Center (SFB) 676 Particles, Strings and the Early Universe.
\section{About subgroups of Sp(1,\textit{n})Sp(1)}
\begin{lemma}[Goursat's theorem]
\label{goursatL} 
Let $\g_1$, $\g_2$ be Lie algebras. There is a one-to-one correspondence between Lie subalgebras $\hh\subset \g_1\oplus \g_2$ and quintuples ${\mathcal Q}(\hh)=(A,A_0,B,B_0,\theta)$, with $A\subset\g_1$ $B\subset\g_2$ Lie subalgebras, $A_0\subset A$, $B_0\subset B$ ideals and $\theta: A/A_0\to B/B_0$ is a Lie algebra isomorphism.
\end{lemma}
\noindent\textit{Proof:} Let $\hh\subset\g_1\oplus\g_2$ be a Lie subalgebra and denote by $\pi_i:\g_1\oplus\g_2\rightarrow \g_i$, $i=1,2$, the natural projections. Set $A:=\pi_1(\hh)\subset \g_1$, $B:=\pi_2(\hh)\subset\g_2$, $A_0:=\ker({\pi_2}_{|\hh})$ and $B_0:=\ker({\pi_1}_{|\hh})$. It is not hard to see that $A_0$ and $B_0$ can be identified with ideals in $A$ and $B$ respectively. Now we can define a map $\tilde{\theta}:A\to B/B_0$ as follows. For $X\in A$ take any $Y\in B$ such that $X+Y\in \hh$ and define $\tilde{\theta}(X):=Y+B_0$. It is easy to check that this map is well defined. Its kernel is $A_0$ so $\tilde{\theta}$ induces a Lie algebra isomorphism $\theta:A/A_0\to B/B_0$. This defines a map $\hh\mapsto{\mathcal Q}(\hh)$.\\
Conversely, a quintuple $Q=(A,A_0,B,B_0,\theta)$ as above defines a Lie subalgebra $\hh={\mathcal G}(Q)\subset\g_1\oplus\g_2$ by setting
$$ \hh:=\lbrace X+Y\in A\oplus B\ | \ \theta(X+A_0)=Y+B_0\rbrace. $$
It is not hard to see that the maps ${\mathcal G}$ and ${\mathcal Q}$ are inverse to each other. $\hfill\Box$

\noindent 
We will use the following two classification results for $\mathbb{H}$-irreducible subgroups of  $\Sp(1,n)$. 
\begin{theorem}[{\cite[Corollary 2.1]{CortesMeinke}}]
\label{ko:irreducibleSp(1,n)}
Let $H\subset \Sp(1,n)$ be a connected and $\HH$-irreducible Lie subgroup. Then $H$ is conjugate to one of the following groups:
\begin{itemize}
	\item[$(i)$] $\SO^0(1,n)$, $\SO^0(1,n)\cdot \U(1)$, $\SO^0(1,n)\cdot \Sp(1)$ if $n\geq 2$,
	\item[$(ii)$] $\SU(1,n)$, $\U(1,n)$,
	\item[$(iii)$] $\Sp(1,n)$,
	\item[$(iv)$] $U^0=\left\{A\in \Sp(1,1)\ |\ A\Phirm=\Phirm A\right\}\cong\Spin^0(1,3)$ with $\Phirm=${\footnotesize$\begin{pmatrix} 0 & -1\\ 1 & 0 \end{pmatrix}$} if $n=1$.
\end{itemize}
\end{theorem} 
\begin{proposition}[{\cite[Proposition 2.4]{CortesMeinke}}]
\label{prop:NormalSubgroupOfIrreducibleGroup}
Let $H\subset\Sp(1,n)$ be an $\HH$-irreducible subgroup. Then one of the following is true.
\begin{itemize}
	\item[$(i)$] $H$ is discrete.
	\item[$(ii)$] $H^0=\U(1)\cdot\mathbbm{1}_{n+1}$ or $H^0=\Sp(1)\cdot\mathbbm{1}_{n+1}$.
	\item[$(iii)$] $H^0$ is $\HH$-irreducible.
		\item[$(iv)$] $n=1$ and $H^0$ is one of the groups $\SO^0(1,1)$, $\SO^0(1,1)\cdot\U(1)$, $\SO^0(1,1)\cdot\Sp(1)$ or $$ S=\left\{\left. e^{ibt}\begin{pmatrix} \cosh(at) & \sinh(at)\\ \sinh(at) & \cosh(at)\end{pmatrix} \right|\ t\in\R \right\}, $$
		for some non-zero real numbers $a,b$.
\end{itemize}
\end{proposition}
\noindent
We denote by $\pi_1:\ssp(1,n)\oplus\ssp(1)\to\ssp(1,n)$ and $\pi_2:\ssp(1,n)\oplus\ssp(1)\to\ssp(1)$ the canonical projections.
\begin{proposition}\label{prop:ListLieAlgebras}
Let $n\geq 2$ and $H\subset\Sp(1,n)\Sp(1)$ be an $\HH$-irreducible closed subgroup. Then the Lie algebra $\hh$ is one of the following:
\begin{itemize}
  \item[$(i)$] $\hh=\hh_0\oplus \cc$ with $\hh_0\in\lbrace \lbrace 0\rbrace,\so(1,n)\rbrace$, $\cc\subset \ssp(1)\cdot\mathbbm{1}_{n+1}\oplus\ssp(1)$ and $\pi_1(\cc)=\ssp(1)\cdot\mathbbm{1}_{n+1}$, $\pi_2(\cc)=\ssp(1)$, $\cc\cap \ssp(1,n)=\lbrace 0\rbrace$, $\cc\cap \ssp(1)=\lbrace 0\rbrace$,
  \item[$(ii)$] $\hh=\hh_0\oplus\cc$ with $\hh_0\in\lbrace \lbrace 0\rbrace, \so(1,n),\su(1,n)\rbrace$, $\cc\subset\uu(1)\cdot\mathbbm{1}_{n+1}\oplus\uu(1)$ and $\pi_1(\cc)=\uu(1)\cdot\mathbbm{1}_{n+1}$, $\pi_2(\cc)=\uu(1)$, $\cc\cap\ssp(1,n)=\lbrace 0\rbrace$, $\cc\cap\ssp(1)=\lbrace 0\rbrace$,
  \item[$(iii)$] $\hh=\hh_0\oplus\cc$ where $\hh_0\subset\ssp(1,n)$ is one of the following Lie algebras
  $$ \ssp(1,n),\ \ \uu(1,n),\ \ \su(1,n),\ \ \so(1,n)\oplus\ssp(1)\cdot\mathbbm{1}_{n+1},\ \ \so(1,n)\oplus\uu(1)\cdot\mathbbm{1}_{n+1},$$ 
  $$\ \ \so(1,n),\ \ \ssp(1)\cdot\mathbbm{1}_{n+1}, \ \ \uu(1)\cdot\mathbbm{1}_{n+1},\ \ \lbrace 0\rbrace, $$
  and $\cc\subset\ssp(1)$ is $\lbrace 0\rbrace, \uu(1)$ or $\ssp(1)$.
\end{itemize}
\end{proposition}
\noindent\textit{Proof:} The idea is to apply Goursat's theorem (Lemma \ref{goursatL}) to $\hh\subset \ssp(1,n)\oplus\ssp(1)$. The Lie subalgebras $A$, $A_0$, $B$ and $B_0$ are given by $\pi_1(\hh)$, $\hh\cap\ssp(1)$, $\pi_2(\hh)$ and $\hh\cap\ssp(1)$. Let $p:\Sp(1,n)\times\Sp(1)\to\Sp(1,n)$ be the natural projection. Notice that $H\subset\Sp(1,n)\Sp(1)$ is $\HH$-irreducible if and only if $p(\hat{H})\subset\Sp(1,n)$ is $\HH$-irreducible, where $\hat{H}$ is the preimage of $H$ under the two-fold covering $\Sp(1,n)\times \Sp(1)\to\Sp(1,n)\Sp(1)$. By Proposition \ref{prop:NormalSubgroupOfIrreducibleGroup} and Theorem \ref{ko:irreducibleSp(1,n)} we know that $p(\hat{H})$ is either discrete or $(p(\hat{H}))^0$ is one of the following subgroups of $\Sp(1,n)$:
$$ \Sp(1,n),\ \  \U(1,n),\ \  \SU(1,n),\ \ \SO^0(1,n)\left( \Sp(1)\cdot\mathbbm{1}_{n+1}\right),\ \ \SO^0(1,n)\left(\U(1)\cdot\mathbbm{1}_{n+1}\right),$$
$$\SO^0(1,n),\ \ \Sp(1)\cdot\mathbbm{1}_{n+1},\ \ \U(1)\cdot\mathbbm{1}_{n+1}. $$
Since $dp=\pi_1$ we immediately obtain all possibilities for $\pi_1(\hh)$. Furthermore $\hh\cap\ssp(1,n)$ is an ideal of the Lie algebra $\pi_1(\hh)$. We can read off from the above list a decomposition of $\pi_1(\hh)$ into ideals, which gives us all possibilities for $\hh\cap\ssp(1,n)$. The resulting list of pairs $(A,A_0)$ is displayed in a table below.\\
On the other side there are only three Lie subalgebras of $\ssp(1)$, namely $\ssp(1)$ itself, $\uu(1)$ and $\lbrace 0\rbrace$. It follows that $\pi_2(\hh)$ is one of these three. Again, $\hh\cap\ssp(1)$ is an ideal of $\pi_2(\hh)$. It follows that the only possibilites for $\hh\cap\ssp(1)$ are the same as for $\pi_2(\hh)$.\\
By Goursat's theorem we have a Lie algebra isomorphism $\theta: A/A_0\to B/B_0$. Since we know all possibilities for $B$ and $B_0$, it follows that $A/A_0$ is isomorphic to $\ssp(1)$, $\uu(1)$ or $\lbrace 0\rbrace$. Therefore we need to consider all possibilities for $A$ and $A_0$, as listed in the following table, and keep only those for which $A/A_0$ is isomorphic to $\ssp(1)$, $\uu(1)$ or $\lbrace 0\rbrace$.
\begin{longtable}{ c | c }
  $A$ & $A_0$ \\
  \hline
  $\ssp(1,n)$ & \begin{tabular}{c} $\ssp(1,n)$\\ $\lbrace 0\rbrace$ \end{tabular}\\
  \hline
  $\su(1,n)\oplus\uu(1)$ & \begin{tabular}{c} $\su(1,n)\oplus\uu(1)$\\ $\su(1,n)$ \\ $\uu(1)$ \\ $\lbrace 0\rbrace$ \end{tabular}\\
  \hline
  $\su(1,n)$ & \begin{tabular}{c} $\su(1,n)$\\ $\lbrace 0\rbrace$ \end{tabular}\\
  \hline
  $\so(1,n)\oplus\ssp(1)$ & \begin{tabular}{c} $\so(1,n)\oplus\ssp(1)$\\ $\so(1,n)$\\ $\ssp(1)$\\ $\lbrace 0\rbrace$ \end{tabular}\\
  \hline
  $\so(1,n)\oplus\uu(1)$ & \begin{tabular}{c} $\so(1,n)\oplus\uu(1)$\\ $\so(1,n)$\\ $\uu(1)$\\ $\lbrace 0\rbrace$ \end{tabular}\\
  \hline
  $\so(1,n)$ & \begin{tabular}{c} $\so(1,n)$\\ $\lbrace 0\rbrace$ \end{tabular}\\
  \hline
  $\ssp(1)$ & \begin{tabular}{c} $\ssp(1)$\\ $\lbrace 0\rbrace$ \end{tabular}\\
  \hline
  $\uu(1)$ & \begin{tabular}{c} $\uu(1)$\\ $\lbrace 0\rbrace$ \end{tabular}\\
  \hline
  $\lbrace 0\rbrace$ & $\lbrace 0\rbrace$
\end{longtable}
\noindent
If $B/B_0\cong\ssp(1)$ then $B=\ssp(1)$ and $B_0=\lbrace 0\rbrace$. The possibilities for $(A,A_0)$ are
$$(\so(1,n)\oplus\ssp(1)\cdot\mathbbm{1}_{n+1}, \so(1,n)) \ \ \textrm{and} \ \ \ (\ssp(1)\cdot\mathbbm{1}_{n+1},\lbrace 0\rbrace).$$ This gives us case $(i)$. Analogously we get the remaining Lie algebras in $(ii)$ and $(iii)$. $\hfill\Box$
\section{Main results}
\subsection{Proof of the main theorem}
\begin{lemma}[{\cite[Lemma 3.1]{CortesMeinke}}]\label{le:InvariantForms}
Let $n\geq 3$ and $\alpha\in \otimes^3V^*$, where $V=\HH^{1,n}$ is considered as real vector space. If $\alpha$ is $\SO^0(1,n)$-invariant, then $\alpha=0$.
\end{lemma}
\begin{remark}\label{re:InvariantForms}
The $\SO^0(1,n)$-invariant elements of $\otimes^3V^*$ are in one-to-one correspondence to the $\SO^0(1,n)$-equivariant bilinear maps from $V\times V$ to $V$. It follows from Lemma $\ref{le:InvariantForms}$ that the corresponding bilinear maps also vanish.
\end{remark}
\noindent\textit{Proof of Theorem \ref{th:MainTheorem}:} Let $\rho: H\to \GL(T_pM)$ be the isotropy representation. We identify $H$ with its image $\rho(H)$. Since $H$ preserves the metric $g$ and the almost quaternionic structure $Q$, we can consider $H$ as a subgroup of $\Sp(1,n)\Sp(1)$.\\
In our first step we consider the covering $G/H^0$ of $M=G/H$ and show that it is a reductive homogeneous space, i.e.\ there exists an $H^0$-invariant subspace $\mathfrak{m}\subset\g$ such that $\g=\hh\oplus\mathfrak{m}$.\\
We apply Proposition \ref{prop:ListLieAlgebras} to $H^0$. The existence of a subspace $\mathfrak{m}$ is clear if $\hh$ is one of the semi-simple Lie algebras in the list. If $\hh$ is one of the abelian Lie algebras contained in $\uu(1)\cdot\mathbbm{1}_{n+1}\oplus\uu(1)$, then the closure of $\Ad(H^0)\subset \GL (\g )$ is compact and hence there exists an $\Ad(H^0)$-invariant subspace $\mathfrak{m}$. The remaining Lie algebras in the list have the form $\hh=\s\oplus \mathfrak{z}$ where $\s$ is semi-simple containing $\so(1,n)$ and $\mathfrak{z}$ is the non-trivial centre. Then $\g$ decomposes into $\g=\s\oplus\mathfrak{z}\oplus\mathfrak{m}$ with respect to the action of $\s$. If we consider the action of $\s$ on $\mathfrak{m}\cong\HH^{1,n}$ as a complex representation, then $\mathfrak{m}$ is either $\C$-irreducible or decomposes into two $\C$-irreducible subrepresentations. Since the elements of $\mathfrak{z}$ commute with $\s$, they preserve the sum of all non-trivial $\s$-submodules, which is precisely $\mathfrak{m}$. Thus we have shown that $G/H^0$ is a reductive homogeneous space.\\
Next we show that $\g=\hh\oplus\mathfrak{m}$ is a symmetric Lie algebra. It is sufficient to show that $\left[ \mathfrak{m},\mathfrak{m}\right]\subset\hh$. We restrict the Lie bracket $\left[\cdot,\cdot\right]$ to $\mathfrak{m}\times\mathfrak{m}$ and denote its projection to $\mathfrak{m}$ by $\beta$. It is an antisymmetric bilinear map which is $\Ad(H)$-equivariant. Since $\mathfrak{m}\cong \HH^{1,n}$, we can consider $\beta$ as an element of $\otimes^3(\HH^{1,n})^*$. It is also $H^{Zar}$-invariant, where $H^{Zar}$ denotes the Zariski closure. Since $H^{Zar}$ is an algebraic group, it has only finitely many connected components, see \cite{Milnor}. Now we show that $(H^{Zar})^0$ is non-compact.\\
Assume that $(H^{Zar})^0$ is compact. Since $H^{Zar}$ has only finitely many connected components it follows that $H^{Zar}$ is compact and therefore contained in a maximal compact subgroup of $\Sp(1,n)\Sp(1)$. Hence, $H^{Zar}$ is conjugate to a subgroup of $(\Sp(1)\times\Sp(n))\Sp(1)$ but this contradicts the $\HH$-irreducibility of $H^{Zar}$. So we have shown that $(H^{Zar})^0$ is non-compact.\\
Now we apply Proposition \ref{prop:ListLieAlgebras} to $H^{Zar}$. Since $H^{Zar}$ is non-compact we see from the list there that $(H^{Zar})^0$ contains $\SO^0(1,n)$. Hence, $\beta$ is $\SO^0(1,n)$-equivariant. Since $n\geq 3$ it follows from Remark \ref{re:InvariantForms} that $\beta$ vanishes. This shows that $\g=\hh\oplus\mathfrak{m}$ is a symmetric Lie algebra and that the universal covering $\tilde{M}=\tilde{G}/\tilde{G}_p$ of $M$ is a symmetric space. The fundamental $4$-form $\Omega$ of $\tilde{M}$ is $\tilde{G}$-invariant and since $\tilde{M}$ is a symmetric space $\Omega$ is parallel. In particular $\Omega$ is closed. It is known that for dimension $\geq 12$ an almost quaternionic Hermitian manifold is quaternionic K\"ahler if $d\Omega=0$, see \cite{Swann}. This shows that $\tilde{M}$ is furthermore a quaternionic K\"ahler manifold. Summarizing, we have shown that $M$ is locally isometric to a quaternionic K\"ahler symmetric space. $\hfill\Box$
\subsection{A class of non-symmetric examples in dimension 12}
\label{counterSec} 
In Theorem \ref{th:MainTheorem} we did not consider the dimension $12$. This is because the arguments used in the proof to show that $M$ is a reductive homogeneous space do not apply in this dimension, although still $\SO^0(1,n)\subset H^{Zar}$ holds. In fact, the proof relies on Lemma \ref{le:InvariantForms} which holds for dimension $4n+4\geq 16$. If $\dim M=12$ then $n=2$ and then there exist non-trivial anti-symmetric bilinear forms $\HH^{1,2}\times\HH^{1,2}\to\HH^{1,2}$ which are invariant under $\SO^0(1,2)$. Therefore in dimension $12$ we cannot be sure if the manifolds are symmetric.\\
In the following we will give a non-symmetric example by specifying a Lie algebra $\g=\hh\oplus\mathfrak{m}$ where $\hh$ is a Lie algebra of the list in Proposition \ref{prop:ListLieAlgebras}. The pair $(\g,\hh)$ defines a simply connected homogeneous space $M=G/H$ where $G$ is a connected and simply connected Lie group with Lie algebra $\g$ and $H$ is the closed connected Lie subgroup of $G$ with Lie algebra $\hh$.\\[3pt]
Let $\hh=\so(1,2)\oplus\cc$ with $\cc=\lbrace (X\cdot\mathbbm{1}_3,X)\in\ssp(1)\cdot\mathbbm{1}_3\oplus\ssp(1)\ |\ X\in\ssp(1) \rbrace$, see Proposition~\ref{prop:ListLieAlgebras}~$(i)$. Then we consider the vector space direct sum $\g :=\hh\oplus\mathfrak{m}$ with $\mathfrak{m}=\HH^{1,2}$ and define a Lie bracket on $\g$ in the following way. For elements $A,B\in\hh$ we take the standard Lie bracket of $\hh$, i.e.\ $\left[A,B \right]=AB-BA$. Then we define $\left[ A,x\right]=-\left[x,A \right]=Ax$ for $A\in\hh$ and $x\in\mathfrak{m}$. Note that, as an $\hh$-module, we can decompose $\mathfrak{m} = \HH^{1,2}=\mathbb{R}^{1,2}\otimes \mathbb{H}=\mathbb{R}^{1,2}\otimes \mathbb{R}^4$, where the action of $\so(1,2)$ is by the defining representation on the first factor and trivial on the second and the action of $\cc \cong \so(3)\subset \so(4)$ is trivial on the first factor and by the standard four-dimensional representation $\mathbb{H} = \mathbb{R} \oplus \mathrm{Im}\, \mathbb{H} = \mathbb{R} \oplus \mathbb{R}^3$ on the second. Finally we have to define the Lie bracket for elements in $\mathfrak{m} = \R^{1,2}\otimes\R^4$.\\
Let $K:\R^{1,2}\to\so(1,2)$ be an isomorphism of Lie algebras where $\R^{1,2}$ is endowed with the Lorentzian cross product, $\iota:\ssp(1)\to\cc$, $X\to X\cdot\mathbbm{1}_3+X$, and let $\eta$ be the standard Lorentz metric on $\R^{1,2}$. Furthermore denote $\left\langle \cdot,\cdot\right\rangle$ the standard inner product on $\R^4$. Let $x=u\otimes p$, $y=v\otimes q\in\R^{1,2}\otimes\R^4$ and write $p=p_0+\vec{p}$, $q=q_0+\vec{q}$, where $p_0, q_0\in \mathbb{R}$ and $\vec{p}, \vec{q}\in \mathrm{Im}\, \mathbb{H} = \mathbb{R}^3$. We set
$$ \left[x,y \right] = \underbrace{\left\langle \vec{p},\vec{q}\right\rangle\cdot K(u\times v)  -\frac12 \eta(u,v)\iota(\vec{p}\times\vec{q})}_{\in\, \hh} + \underbrace{u\times v(p_0q_0-\left\langle \vec{p},\vec{q}\right\rangle )}_{\in\, \mathbb{R}^{1,2}\subset\, \mathbb{H}^{1,2}\, =\, \mathfrak{m}}, $$
where $\vec{p}\times\vec{q}$ is the Euclidian cross product in $\mathrm{Im}\, \mathbb{H}=\ssp(1)$ and $u\times v$ the Lorentzian cross product in $\mathbb{R}^{1,2}$. This extends the partially defined bracket to an anti-symmetric bilinear map $\left[\cdot,\cdot\right]:\g\times\g\to\g$, which satisfies the Jacobi-identity. Hence $\g$ becomes a Lie algebra. We claim that $(\g ,\hh )$ is not a symmetric pair. In fact, every $\hh$-invariant complement $\mathfrak{m}'$ of $\hh$ in $\g$ contains $\mathbb{R}^{1,2}\otimes \mathbb{R}^3$ (there is no other equivalent $\hh$-submodule in $\g$) and thus we see from the formula for the bracket that $\left[\mathfrak{m}',\mathfrak{m}' \right]\nsubseteq \hh$.\\

For a general classification of the homogeneous spaces with $\hh=\so(1,2)\oplus\cc$ we need to classify all the Lie algebra structures on the vector $\g=\hh\oplus\R^{1,2}\otimes\R^4$ such that the Lie bracket restricts to the Lie bracket of $\hh$ and to the given representation of $\hh$ on $\R^{1,2}\otimes\R^4$. For this one has to describe all the $\hh$-invariant tensors of $\mathrm{\Lambda}^2\mathfrak{m}^*\otimes\g\cong \mathrm{\Lambda}^2\mathfrak{m}^*\otimes\hh\oplus\mathrm{\Lambda}^2\mathfrak{m}^*\otimes\mathfrak{m}$ which satisfy the Jacobi-identity. With the above notation, these bilinear maps have the following form
\begin{eqnarray*}
  \left[x,y \right] &=& \left( a\cdot p_0q_0+b\left\langle \vec{p},\vec{q}\right\rangle \right)\cdot K(u\times v) + \eta(u,v)\left( c\cdot\iota(\vec{p}\times\vec{q})+d\left( p_0\vec{q}-q_0\vec{p}\right) \right)\\
  & & +u\times v\cdot \left( a_1\cdot p_0q_0+a_2\cdot\left\langle \vec{p},\vec{q} \right\rangle + \frac{a_3}{2}\left( p_0\vec{q}+q_0\vec{p}\right) \right),
\end{eqnarray*}
where $a,b,c,d,a_1,a_2,a_3\in\R$. The bracket satisfies the Jacobi-identity if and only if the following equations hold
\begin{eqnarray}
  0&=& d,\nonumber\\
  0&=&a+\frac{a_1a_3}{2}-\frac{a_3^2}{4},\\
  0&=&b+2c+\frac{a_2a_3}{2},\\
  0&=&b+a_1a_2-\frac{a_2a_3}{2},\\
  0&=&-\frac{ba_3}{2}+aa_2.
\end{eqnarray}
Summarizing we obtain the following proposition. 
\begin{proposition} \label{classProp} Every solution $(a,b,c,a_1,a_2,a_3)$ of the 
quadratic system (1)-(4) defines a connected and simply connected 
homogeneous almost quaternionic pseudo-\-Hermitian
manifold $G/H$ with isotropy algebra $\hh=\so(1,2)\oplus\so(3)\subset \so(1,2)\oplus\so(4)\subset \mathfrak{gl}(\mathbb{R}^{1,2}\otimes \mathbb{R}^4)\cong \mathfrak{gl}(12,\mathbb{R})$. Conversely, every 
such homogeneous space arises by this construction. 
\end{proposition}
\noindent The above example corresponds to $a=0$, $b=1$, $c=-\frac{1}{2}$, $d=0$, $a_1=1$, $a_2=-1$ and $a_3=0$.

\end{document}